\documentclass[11pt,twoside]{article}
 \usepackage{amssymb,amsfonts}
 \usepackage{amsmath}

\newtheorem{theorem}{Theorem}[section]
\newtheorem{corollary}[theorem]{Corollary}
\newtheorem{lemma}[theorem]{Lemma}
\newtheorem{remark}[theorem]{Remark}

\newtheorem{example}[theorem]{Example}

\numberwithin{equation}{section}
\newcommand{\field}[1]{\mathbb{#1}}

\newcommand{\R}{\field{R}}

\begin{document}

\title{Global attractor in competitive Lotka-Volterra systems}

\author{Zhanyuan Hou}

\date{}
\maketitle
\begin{abstract}
For autonomous Lotka-Volterra systems of differential equations
modelling the dynamics of $n$ competing species, new criteria are
established for the existence of a single point global attractor.
Under the conditions of these criteria, some of the species will
survive and stabilise at a steady state whereas the others, if
any, will die out.
\end{abstract}

\textbf{Keywords} Lotka-Volterra, global attractor, competitive
systems

\textbf{MSC(2000)} 34D45, 92D25

\section{Introduction}

In this paper, we are concerned with the asymptotic behaviour of
solutions of the autonomous Lotka-Volterra system
\begin{equation}
x'_i = b_ix_i(1-\alpha_ix), \quad i\in I_n, \label{1.1}
\end{equation}
where $x'_i = dx_i/dt$, $b_i >0$, $I_m =\{1, 2, \ldots, m\}$ for
any positive integer $m$, $\alpha_i = (a_{i1}, a_{i2}, \ldots,
a_{in})$ with $a_{ii}>0$ and $a_{ij}\geq 0$, and $x = (x_1, x_2,
\ldots, x_n)^T \in {\R}^n$. Since the behaviour of $(\ref{1.1})$
is clear when $n=2$, we assume $n>2$ though some of our results
given later are applicable to two-dimensional systems. Let
${\R}^n_+ = \{x\in{\R}^n: x\geq 0\}$, where by $y\geq x$ in
${\R}^n$ we mean $y_i\geq x_i$ for all $i\in I_n$. We also write
$y >x$ if $y\geq x$ but $y\not= x$. We shall restrict the
solutions of $(\ref{1.1})$ to ${\R}^n_+$ in view of $(\ref{1.1})$
modelling the dynamics of $n$ competitive species, $x_i(t)$
denoting the population size of the $i$th species at time $t$. We
say that an equilibrium $x^*\in {\R}^n_+$ of $(\ref{1.1})$ is a
global attractor if every solution of $(\ref{1.1})$ in
int${\R}^n_+ = \{x\in {\R}^n: x_i>0, i\in I_n\}$, i.e. with
$x(0)\in \textrm{int}{\R}^n_+$, satisfies $\lim_{t\to\infty}x(t) =
x^*$.

Since the independent introduction of $(\ref{1.1})$ as a
population model by Lotka and Volterra in the 1920s, a large
number of researchers have been attracted to the theoretical
investigation of $(\ref{1.1})$ and its variations (see, for
example, \cite{Ze1} and the references therein for a brief account
of the development). One of the problems arising from
$(\ref{1.1})$ as a model is to find suitable conditions under
which some particular (or all) species in the competitive
community will survive and stabilise at a steady state whereas the
rest, if any, will die out eventually. Mathematically, we need
suitable conditions which ensure the existence of an equilibrium
$x^*\in{\R}^n_+$ that is a global attractor. It is known that
$(\ref{1.1})$ has an equilibrium $x^*\in \textrm{int}{\R}^n_+$
that is a global attractor if
\begin{equation}
\frac{a_{i1}}{a_{11}} + \frac{a_{i2}}{a_{22}} + \cdots +
\frac{a_{in}}{a_{nn}} <2 \label{1.2}
\end{equation}
holds for all $i\in I_n$ (see \cite{Ka}, \cite{GoAh} or \cite[p.
294--297]{Go}). For $(\ref{1.1})$ to have the same property as
above, Hou \cite{Ho1} finds a more general condition
\begin{equation}
\max\biggl\{0, \frac{a_{ij}}{a_{jj}}\biggl(1-\sum_{k\not = i, j}
\frac{a_{jk}}{a_{kk}}\biggr)\biggr\} < 1 - \sum_{k\not = i, j}
\frac{a_{ik}}{a_{kk}} \label{1.3}
\end{equation}
for all $i, j\in I_n$ with $i\not = j$, which incorporates
$(\ref{1.2})$ for all $i\in I_n$ as a special case. For
$(\ref{1.1})$ to have an equilibrium $x^* \in
\partial{\R}^n_+ = {\R}^n_+\setminus \textrm{int}{\R}^n_+$ that
is a global attractor, Zeeman \cite{Ze2} shows that
\begin{equation}
(i-j)(a_{ij}-a_{jj}) >0 \label{1.4}
\end{equation}
for $i, j\in I_n$ with $i\not = j$ guarantees the survival of only
the first species, and Ahmad and Lazer \cite{AhLa} obtain a
criterion for the extinction of only one species. Li, Yu and Zeng
\cite{LiYuZe} extend this criterion to the extinction of possibly
more than one species. Combining $(\ref{1.3})$ with a condition
similar to that of \cite{MoZe} for nonautonomous systems, Hou
\cite{Ho1} shows the global attraction of an equilibrium $x^*
\in\partial {\R}^n_+$ with the survival of some species and the
extinction of the others.

The main purpose of this paper is to establish new criteria for
the existence of a single point global attractor. Using the idea
of the ``ultimate contracting cells'' method, which was summarised
by Hou \cite{Ho2} or \cite{Ho3}, rather than using Lyapunov
functions, we shall see that the new conditions are geometrically
simple and include $(\ref{1.2})$, $(\ref{1.3})$ and those given in
\cite{AhLa} as particular instances.

The presentation of this paper is as follows. The statements of
the main results under geometric conditions, together with some
examples, will be given in \S2. Then in \S3 these geometric
conditions will be translated into algebraic expressions. In \S4
the main results will be restated under equivalent algebraic
conditions. Finally, the proofs of the main theorems will be
completed in \S5.

\section{Main results under geometric conditions}

We first explain some terms and notation. For any $u,
v\in{\R}^n_+$ with $v\geq u$, a \textit{cell} $[u, v]$ is defined
to be the set of $x\in{\R}^n_+$ satisfying $v\geq x\geq u$. For
any $J\subset I_n$, $u^J\in {\R}^n_+$ is defined by $u^J_i = u_i$
if $i\in J$ and $u^J_i = 0$ otherwise. Thus, with this notation,
$u^{I_n} = u$, $u^{\emptyset} = 0$, and $u^{I_n\setminus\{k\}}$
for any $k\in I_n$ is the projection of $u$ onto the coordinate
plane $\pi_k$, which is defined by $(\ref{2.2})$. If $u
> 0$, the set $\gamma = \{x\in{\R}^n_+: u^Tx=1\}$ is an
$(n-1)$-dimensional \textit{plane} in ${\R}^n_+$. A point
$x\in{\R}^n_+$ is said to be \textit{below} (\textit{on} or
\textit{above}) $\gamma$ if $u^Tx <1$ ($u^Tx =1$ or $u^Tx >1$). A
nonempty set $S\subset {\R}^n_+$ is said to be \textit{below}
(\textit{on} or \textit{above}) $\gamma$ if every point $x$ in $S$
is so. Now for $(\ref{1.1})$ we let
\begin{eqnarray}
Y &=& (a_{11}^{-1}, a_{22}^{-1}, \ldots, a_{nn}^{-1})^T, \label{2.1} \\
\pi_i &=& \{x\in{\R}^n_+: x_i=0\}, \quad i\in I_n, \label{2.2}\\
\gamma_i &=&  \{x\in{\R}^n_+: \alpha_ix=1\}, \quad i\in I_n,
\label{2.3}
\end{eqnarray}
and view $\cap_{i\in\emptyset}\pi_i$ as ${\R}^n_+$. Then, for any
$J\subset I_n$, $\cap_{i \in J}\pi_i$ is obviously invariant under
$(\ref{1.1})$. Thus, as $\gamma_i\cap(\cap_{j \in
I_n\setminus\{i\}}\pi_j) = \{Y^{\{i\}}\}$, $Y^{\{i\}}$ is an
equilibrium of $(\ref{1.1})$ on the $x_i$-axis. Let $U\in [0, Y]$
be defined as follows: for each $i\in I_n$, if $\gamma_i$ is below
$\gamma_j$ for every $j\in I_n\setminus\{i\}$ then $U_i = 0$; if
$Y^{\{i\}}$ is on or above some $\gamma_j$ $(j\not = i)$ then $U_i
= Y_i = a^{-1}_{ii}$; otherwise $U_i$ is the maximum value of the
$i$th coordinate components of the intersection points of
$\gamma_i$ with any other $\gamma_j$.

\vskip 4 mm \textbf{CONDITION $(C_k)$.} Either
$U^{I_n\setminus\{k\}}$ is below $\gamma_k$ or $\gamma_k\cap[0,
U^{I_n\setminus\{k\}}]$ is above $\gamma_j$ for every $j\in
I_n\setminus\{k\}$. Alternatively, either $U^{I_n\setminus\{k\}}$
is below $\gamma_k$ or the set $\gamma_j\cap[0,
U^{I_n\setminus\{k\}}]$ is below $\gamma_k$ for every $j\in
I_n\setminus\{k\}$.

\vskip 4 mm The analytical definition of $U$, as well as its
simplification, and algebraic realisation of condition $(C_k)$ are
left to the next section. We shall see in \S5 that $U$ is an
ultimate upper bound for every solution $x(t)$ of $(\ref{1.1})$ in
$\textrm{int}{\R}^n_+$, i.e. $\overline{\lim}_{t\to\infty}x_i(t)
\leq U_i$ for all $i\in I_n$. Moreover, condition $(C_k)$ will
guarantee the existence of some $\delta > 0$ such that every
solution of $(\ref{1.1})$ in int${\R}^n_+$ satisfies
$\underline{\lim}_{t\to\infty}x_k(t) \geq \delta$. Most
importantly, the following theorem holds.

\begin{theorem} \label{Th2.1}
Assume that $(C_k)$ holds for all $k\in I_n$. Then $(\ref{1.1})$
has a unique equilibrium $x^*\in \textup{int}{\R}^n_+$ that is a
global attractor. Moreover, the inequality $x^*\leq U$ holds.
\end{theorem}

\begin{remark} \label{Re2.2}
\textup{It is shown in \cite{Ho1} that $(\ref{1.2})$ for all $i
\in I_n$ implies $(\ref{1.3})$ for all $i, j\in I_n$ with $i\not
=j$. It is also shown that $(\ref{1.3})$ for all $i, j\in I_n$
with $i\not =j$ is equivalent to that $\gamma_j\cap[0,
Y^{I_n\setminus \{k\}}]$ is below $\gamma_k$ for all $j, k\in I_n$
with $j\not= k$. Since $U\leq Y$, $(\ref{1.3})$ for all $i, j\in
I_n$ with $i\not= j$ implies condition $(C_k)$ for all $k\in I_n$.
Therefore, theorem \ref{Th2.1} incorporates the known results
using $(\ref{1.2})$ and $(\ref{1.3})$ as special instances. The
example below satisfies $(C_k)$ for all $k\in I_n$ but does not
satisfy $(\ref{1.3})$ for some $i, j\in I_n$. Hence, theorem
\ref{Th2.1} applies to a broader class of systems than the known
results do.}
\end{remark}

\begin{example}\label{Ex2.3}
\textup{Consider system $(\ref{1.1})$ with $n=3$ and
\begin{equation} \label{2.4}
\left(\begin{array}{ccc}
a_{11} & a_{12} & a_{13} \\
a_{21} & a_{22} & a_{23} \\
a_{31} & a_{32} & a_{33} \end{array}\right) =
\left(\begin{array}{ccc}
2.1 & 1.9 & 1.9 \\
1 & 3 & 1 \\
1 & 1 & 3 \end{array}\right).
\end{equation}
Then $a_{ii} > a_{ji}$ holds for $i, j\in I_3$ with $i\not= j$ so
$Y^{\{i\}}$ is below $\gamma_j$ for all $i, j\in I_3$ with
$i\not=j$. Since $\gamma_1\cap\gamma_2\cap\pi_3 = \{(0.25, 0.25,
0)^T\}$, $\gamma_1\cap\pi_2\cap\gamma_3 = \{(0.25, 0, 0.25)^T\}$,
$\pi_1\cap\gamma_2\cap\gamma_3 = \{(0, 0.25, 0.25)^T\}$ and
$(\cap_{j \in I_3\setminus\{i\}}\gamma_j)\cap\pi_i$ is below
$\gamma_i$ for $i\in I_3$, by definition we have $U = (0.25, 0.25,
0.25)^T$. Note that $\alpha_iU^{I_3\setminus\{i\}} < 1$, so
($C_i$) holds, for all $i\in I_3$ (see Figure 2.1). By theorem
\ref{Th2.1}, $\gamma_1\cap\gamma_2\cap\gamma_3 = \{x^*\}$ with
$U_i\geq x^*_i>0$ for $i\in I_3$ and $x^*$ is a global attractor.
Indeed, $x^* = (\frac{1}{23}, \frac{11}{46}, \frac{11}{46})^T$. On
the other hand, however, for $Y = (\frac{10}{21}, \frac{1}{3},
\frac{1}{3})^T$, we have $(0, \frac{1}{3}, \frac{11}{57})^T \in
\gamma_1 \cap[0, Y^{\{2, 3\}}]$ but $(0, \frac{1}{3},
\frac{11}{57})^T$ is below $\gamma_3$. So $(\ref{1.3})$ does not
hold for some $i, j\in I_3$.}

\begin{figure}
\begin{picture}(350, 140)(60, 20)
\thicklines \put(80, 40){\vector(1, 0){110}} \put(80,
40){\vector(0, 1){100}} \put(80, 40){\dashbox{1}(30, 30)}
\put(110, 40){\line(-1, 3){30}} \put(80, 70){\line(3, -1){90}}
\put(127, 40){\line(-1, 1){47}} \put(80, 40){\dashbox{1}(22.5,
22.5)}
 \put(72, 32){$O$} \put(68,
56){$U_3$} \put(75, 145){$x_3$} \put(88, 80){$\gamma_1$} \put(90,
110){$\gamma_2$} \put(110, 70){$Y^{\{2, 3\}}$} \put(155,
125){$\pi_1$} \put(160, 50){$\gamma_3$} \put(98, 30){$U_2$}
\put(195, 35){$x_2$} \put(240, 40){\vector(1, 0){110}} \put(240,
40){\vector(0, 1){100}} \put(285, 40){\line(-1, 1){45}} \put(240,
70){\line(3, -1){90}} \put(330, 40){\line(-1, 1){90}} \put(240,
40){\dashbox{1}(22.5, 22.5)} \put(232, 32){$O$} \put(228,
56){$U_3$} \put(235, 145){$x_3$} \put(250, 77){$\gamma_1$}
\put(260, 115){$\gamma_2$} \put(285, 55){$\gamma_3$} \put(315,
125){$\pi_2$} \put(257,
 30){$U_1$} \put(355, 35){$x_1$}
\put(75, 18){Figure 2.1.
 $U^{I_3\setminus\{i\}}$ is below $\gamma_i$ for $i\in I_3$, but
 $\gamma_3\cap[0, Y^{\{2, 3\}}]$ does not below $\gamma_1$.}
\end{picture}
\end{figure}
\end{example}

\begin{corollary} \label{Co2.4} Assume that $(C_k)$
holds for all $k\in I_n$. Then, for each nonempty set $J\subset
I_n$, $(\cap_{j\in J}\gamma_j)\cap(\cap_{\ell\in I_n\setminus
J}\pi_{\ell}) = \{u_J\}$ with $U_j \geq (u_J)_j >0$ for $j\in J$
such that every solution of $(\ref{1.1})$ with $x_j(0) >0$ for
$j\in J$ and $x_{\ell}(0) =0$ for $\ell\in I_n\setminus J$
satisfies $\lim_{t\to\infty}x(t) = u_J$. Moreover, $(\ref{1.1})$
has $2^n$ equilibria in ${\R}^n_+$.
\end{corollary}

\textit{Proof.} When $J = \{i\}$ for some $i\in I_n$,
$(\ref{1.1})$ on $\cap_{k\in I_n\setminus\{i\}}\pi_k$ is
equivalent to $x'_i = b_ix_i(1-a_{ii}x_i)$ so $u_J = Y^{\{i\}}$
meets the requirement. Note that theorem \ref{Th2.1} can be
applied to $(\ref{1.1})$ with $n>1$. Then, if $|J|>1$, the
conclusion about the required $u_J$ follows from the application
of theorem \ref{Th2.1} to the $|J|$-dimensional system
$(\ref{1.1})$ on $\cap_{\ell \in I_n\setminus J}\pi_{\ell}$.
Clearly, $u_J$ is the unique equilibrium of $(\ref{1.1})$ in
$\cap_{\ell \in I_n\setminus J}\pi_{\ell}$ satisfying $(u_J)_j >0$
for $j\in J$. Since $I_n$ has $2^n-1$ nonempty subsets and the
origin is also an equilibrium, $(\ref{1.1})$ has $2^n$ equilibria
in ${\R}^n_+$.

We have mentioned earlier that condition $(C_k)$, for any
particular $k\in I_n$, will secure the survival of the $k$th
species. Bearing this in mind and applying condition $(C_k)$ to
some of the indices (not necessarily all), we have the following.

\begin{theorem} \label{Th2.5}
Assume that condition $(C_k)$ holds for all $k\in J\subset I_n$
$(J\not=\emptyset)$. Then $(\cap_{j\in J}\gamma_j)\cap(\cap_{k\in
I_n\setminus J}\pi_k) = \{x^*\}$ with $0 <x^*_j \leq U_j$ for
$j\in J$. In addition, if $x^*$ is on $\gamma_k$ for every $k\in
I_n\setminus J$, then $x^*$ is a global attractor. Moreover, the
conclusion still holds if ``$x^*$ is on $\gamma_k$'' is replaced
by ``$x^*$ is on or above $\gamma_k$'' when $U$ is replaced by
$Y$.
\end{theorem}

\begin{remark}
\label{Re2.6} \textup{Theorem \ref{Th2.1} is the extreme case of
theorem \ref{Th2.5} when $J = I_n$. When $J=I_{n-1}$, Ahmad and
Lazer \cite{AhLa} obtained the same conclusion under the
assumptions that $(i)$ $(\ref{1.2})$ holds for all $i\in I_{n-1}$,
$(ii)$ $\alpha_nx^*
>1$ and $(iii)$ $\det(a_{ij})_{n\times n} >0$. Later in
\cite{AhLa2, AhLa3} it was pointed out that the above condition
$(iii)$ is redundant. From remark \ref{Re2.2} we know that, for
any $i\in I_n$, $(\ref{1.2})$ implies $(\ref{1.3})$ for all
$j\not=i$, which ensures condition $(C_i)$ with the replacement of
$U$ by $Y$. Thus, the result in \cite{AhLa} under the assumptions
$(i)$ and $(ii)$ is covered by theorem~\ref{Th2.5} as a special
case.}
\end{remark}

\begin{example}
\textup{Consider system $(\ref{1.1})$ with $n=4$ and
\[ \left(\begin{array}{cccc}
a_{11} & a_{12} & a_{13} & a_{14} \\
a_{21} & a_{22} & a_{23} & a_{24} \\
a_{31} & a_{32} & a_{33} & a_{34} \\
a_{41} & a_{42} & a_{43} & a_{44} \end{array}\right) =
\left(\begin{array}{cccc}
5/2 & 3/2 & 1 & 3/2 \\
1 & 3 & 1 & 1 \\
2 & 2 & 2 & 2 \\
2 & 2 & 1 & 3  \end{array}\right). \]
 Then $Y = (\frac{2}{5},
\frac{1}{3}, \frac{1}{2}, \frac{1}{3})^T$. Clearly, $Y^{\{i\}}$ is
below $\gamma_j$ for every $j\in I_4 \setminus\{i\}$ as $a_{ii} >
a_{ji}$. To calculate $U_1$, we find all the possible intersection
points of $\gamma_1$ with $\gamma_j$ ($j\not=1$) on
$\pi_2\cap\pi_3$, $\pi_2\cap\pi_4$ and $\pi_3\cap\pi_4$ and find
the maximum of the $x_1$ values of these points. On
$\pi_3\cap\pi_4$, $\gamma_1$ intersects $\gamma_j$ at
$(\frac{1}{4}, \frac{1}{4}, 0, 0)^T$ for all $j\in
I_4\setminus\{1\}$. On $\pi_2\cap\pi_4$, $\gamma_1$ intersects
$\gamma_3$ at $(\frac{1}{3}, 0, \frac{1}{6}, 0)^T$; $\gamma_1$
intersects $\gamma_2$ and $\gamma_4$ on $\pi_1$ so these are
ignored. On $\pi_2\cap\pi_3$, $\gamma_1$ intersects $\gamma_3$ at
$(\frac{1}{4}, 0, 0, \frac{1}{4})^T$ and $\gamma_4$ at
$(\frac{1}{3}, 0, 0, \frac{1}{9})^T$; $\gamma_1$ is below
$\gamma_2$. Then $U_1 = \max\{\frac{1}{4}, \frac{1}{3}\} =
\frac{1}{3}$. We then obtain $U_2 = \frac{2}{7}$, $U_3 = U_4 =
\frac{1}{4}$ in the same manner. Condition $(C_2)$ holds since
$\alpha_2U^{I_4\setminus\{2\}} = U_1+U_3 +U_4 =\frac{5}{6} <1$.
The set $[0, U^{I_4\setminus\{1\}}]\cap\gamma_1$ is a triangle
with vertices $P_1 = (0, \frac{2}{7}, \frac{1}{4},
\frac{3}{14})^T$, $P_2 = (0, \frac{2}{7}, \frac{11}{56},
\frac{1}{4})^T$ and $P_3 = (0, \frac{1}{4}, \frac{1}{4},
\frac{1}{4})^T$. Since $\alpha_jP_i >1$ for all $i\in I_3$ and
$j\in I_4\setminus\{1\}$, $[0, U^{I_4\setminus\{1\}}]\cap\gamma_1$
is above $\gamma_j$ for every $j\in I_4\setminus\{1\}$. Thus,
condition $(C_1)$ also holds. Then
$\gamma_1\cap\gamma_2\cap\pi_3\cap\pi_4 = \{x^*\}$ with $x^* =
(\frac{1}{4}, \frac{1}{4}, 0, 0)^T$. Since $\alpha_3x^*=1$ and
$\alpha_4x^*=1$, by theorem \ref{Th2.5}, $x^*$ is a global
attractor.}
\end{example}

Another extreme case of theorem \ref{Th2.5} when $J = \{k\}$ is
simplified to the corollary below.

\begin{corollary} \label{Co2.8}
Assume that condition $(C_k)$ and $a_{kk} = a_{jk}$ hold for some
$k\in I_n$ and all $j\in I_n\setminus\{k\}$. Then $Y^{\{k\}}$ is a
global attractor. The conclusion still holds if $U$ is replaced by
$Y$ and $a_{kk} = a_{jk}$ is replaced by $a_{kk} \leq a_{jk}$.
\end{corollary}

\begin{remark}\label{Re2.9}
\textup{The conditions of corollary \ref{Co2.8} with $k=1$ and
$(\ref{1.4})$ for $i, j\in I_n$ with $i\not= j$ are mutually
exclusive.}
\end{remark}

If $U_{\ell} = 0$ for some $\ell\in I_n$, from the definition of
$U$ we make the following observations:
\begin{itemize}
\item[$(i)$] For each $j\in I_n\setminus\{\ell\}$, $Y^{\{\ell\}}$
is below $\gamma_j$ and every $y$ in $\gamma_{\ell}$ is on or
below $\gamma_j$. Hence, $\gamma_j\cap[0,
U^{I_n\setminus\{\ell\}}]$ cannot be below $\gamma_{\ell}$ so
condition $(C_{\ell})$ does not hold. \item[$(ii)$] For $j\in
I_n\setminus\{\ell\}$, $Y^{\{j\}}\in\gamma_j$ so, by $(i)$,
$Y^{\{j\}}$ cannot be below $\gamma_{\ell}$. Thus,
$U^{I_n\setminus\{\ell\}} = Y^{I_n\setminus\{\ell\}}$.
\item[$(iii)$] From lemma~\ref{Le5.4} given in \S5 we shall see
that every solution of $(\ref{1.1})$ in int${\R}^n_+$ satisfies
$\lim_{t\to\infty}x_{\ell}(t) = 0$.
\end{itemize}
Based on these observations we can redefine
$U^{I_n\setminus\{\ell\}}$ by the same definition as before except
the substitution of $I_n\setminus\{\ell\}$ for $I_n$. Repeat this
process until no more new zero components occur.

\begin{theorem} \label{Th2.10}
Let $J = \{j_1, j_2, \ldots, j_m\}\subset I_n$ with $1\leq m<n$
and let $V$ be defined by $V^J=0$ and $V^{I_n\setminus J}\in [0,
Y^{I_n\setminus J}]$ given by the definition of $U$ with the
replacement of $I_n$ by $I_n\setminus J$. Assume that the
following hold:
\begin{itemize}
\item[$(i)$] For each $\ell\in I_m$, $Y^{\{j_{\ell}\}}$ is below
$\gamma_k$ and every $y$ in
$\gamma_{j_{\ell}}\cap(\pi_{j_1}\cap\cdots\cap\pi_{j_{\ell-1}})$
is on or below $\gamma_k$ for all $k\in I_n\setminus\{j_1, j_2,
\ldots, j_{\ell}\}$. \item[$(ii)$] For each $k\in I_n\setminus J$,
either $V^{I_n\setminus\{k\}}$ is below $\gamma_k$ or
$\gamma_k\cap[0, V^{I_n\setminus\{k\}}]$ is above $\gamma_j$ for
all $j\in I_n\setminus(J\cup\{k\})$.
\end{itemize}
Then $(\cap_{j\in J}\pi_j)\cap(\cap_{k\in I_n\setminus J}\gamma_k)
= \{x^*\}$ with $0<x^*_k\leq V_k$ for $k\in I_n\setminus J$ and
$x^*$ is a global attractor.
\end{theorem}

\begin{remark}\label{Re2.11}
\textup{Suppose the conditions of theorem \ref{Th2.10} are met.
Then, applying theorem \ref{Th2.10} to $(\ref{1.1})$ on
$\cap_{\ell\in S}\pi_{\ell}$ for any subset $S\subset J$, we see
that every solution of $(\ref{1.1})$ with $x_k(0)>0$ for $k\in
I_n\setminus J$ and $x_j(0)\geq 0$ for $j\in J$ satisfies
$\lim_{t\to\infty}x(t) = x^*$. This remark also applies to
corollary \ref{Co2.13}.}
\end{remark}

\begin{example}
\textup{Consider system $(\ref{1.1})$ with $n=5$ and
\[
\left(\begin{array}{ccccc}
a_{11} & a_{12} & a_{13} & a_{14} & a_{15} \\
a_{21} & a_{22} & a_{23} & a_{24} & a_{25} \\
a_{31} & a_{32} & a_{33} & a_{34} & a_{35} \\
a_{41} & a_{42} & a_{43} & a_{44} & a_{45} \\
a_{51} & a_{52} & a_{53} & a_{54} & a_{55} \end{array}\right) =
\left(\begin{array}{ccccc}
2.1 & 1.9 & 1.9 & 0 & 1 \\
1 & 3 & 1 & 2 & 0 \\
1 & 1 & 3 & 1 & 1 \\
3 & 3 & 3 & 2.1 & 0 \\
4 & 3 & 3 & 2.5 & 2  \end{array}\right).
\]
Since $a_{5j} \geq a_{ij}$ and $a_{55} > a_{i5}$ for all $i, j\in
I_4$, $Y^{\{5\}}$ is below $\gamma_j$ and each $y\in \gamma_5$ is
on or below $\gamma_j$ for all $j\in I_4$. Also, $a_{4j}\geq
a_{ij}$ and $a_{44} > a_{i4}$ for all $i, j\in I_3$, so
$Y^{\{4\}}$ is below $\gamma_j$ and each $y\in \gamma_4\cap\pi_5$
is on or below $\gamma_j$ for all $j\in I_3$. Now with $J = \{j_1,
j_2\}$, $j_1 = 5$ and $j_2=4$, condition $(i)$ of theorem
\ref{Th2.10} is met. Since the entries $a_{ij}$ for $i, j\in I_3$
are identical with those given in $(\ref{2.4})$, from example
\ref{Ex2.3} we obtain $V=(\frac{1}{4}, \frac{1}{4}, \frac{1}{4},
0, 0)^T$ and $\alpha_iV^{I_5\setminus\{i\}} <1$ for all $i\in
I_5\setminus J = I_3$. By theorem \ref{Th2.10}, the equilibrium
$x^* = (\frac{1}{23}, \frac{11}{46}, \frac{11}{46}, 0, 0)^T$ is a
global attractor.}
\end{example}

When $n = m+1$ with $I_n = J\cup\{k\}$, $V^{I_n\setminus\{k\}} =
0$ so assumption $(ii)$ is automatically met and theorem
\ref{Th2.10} has the following simplified form.

\begin{corollary} \label{Co2.13}
Assume that, for each $\ell\in I_{n-1}$ and every $j\in
I_n\setminus I_{\ell}$, $Y^{\{\ell\}}$ is below $\gamma_j$ and
every $y$ in $\gamma_{\ell}\cap(\pi_1\cap\cdots\cap\pi_{\ell-1})$
is on or below $\gamma_j$. Then $Y^{\{n\}}$ is a global attractor.
\end{corollary}

\begin{remark} \label{Re2.14}
\textup{The condition of corollary \ref{Co2.13} is equivalent to
\begin{equation} \label{2.5}
a_{j\ell}<a_{\ell\ell}, \quad a_{jj}\leq a_{j-1 j} \leq \cdots
\leq a_{1 j}, \quad 1\leq \ell <j\leq n.
\end{equation}
Note that $(\ref{1.4})$ and $(\ref{2.5})$ are mutually exclusive.
Thus, neither of corollary \ref{Co2.13} and the result given by
Zeeman \cite{Ze2} is covered by the other.}
\end{remark}

Finally, we make a brief remark about the conditions of the main
result given in \cite{LiYuZe}. It is assumed that, for some $r<n$,
$(\ref{1.2})$ with the replacement of $n$ by $r$ holds for all
$i\in I_r$ This implies the existence of an equilibrium $x^*\in
\cap_{j\in I_n\setminus I_r}\pi_j$ with $x^*_i>0$ for $i\in I_r$.
It is also assumed that $x^*$ is above $\gamma_j$ for every $j\in
I_n\setminus I_r$. Then, with another condition requiring the
existence of $n-r$ constants satisfying an inequality, it is shown
that $x^*$ is a global attractor. These conditions and those of
theorems \ref{Th2.5} and \ref{Th2.10} are mutually independent
though they have some overlap.

\section{Analytical definition of $U$ and algebraic
expression of $(C_k)$}

The definition of $U$ given in \S2 can be formulated as follows.
For each $i\in I_n$, let
\begin{equation} \label{3.1}
U_i = \left\{\begin{array}{l} a^{-1}_{ii} \quad \textrm{if} \;
a_{ii}\leq a_{ji} \; \textrm{or} \;
a_{ij}=0 \; \textrm{for some} \; j\in I_n\setminus\{i\}, \\
0 \qquad \textrm{if} \; a_{ji} <a_{ii}, a_{jk}\leq a_{ik}\;
\textrm{for all}\; j, k\in I_n\setminus\{i\}, \\
\max\left\{\frac{a_{kj}-a_{ij}}{a_{ii}a_{kj}-a_{ij}a_{ki}}: j,
k\in I_n\setminus\{i\}, a_{kj} > a_{ij}\right\} \quad
\textrm{otherwise}.
\end{array}\right.
\end{equation}

\begin{remark} \label{Re3.1}
\textup{If condition $(C_k)$ holds for some $k\in I_n$, then every
equilibrium of $(\ref{1.1})$ on $\pi_k$ is below $\gamma_k$.
Indeed, $u=0$ is obviously below $\gamma_k$. Let $u\in \pi_k$ with
$u\not=0$ be an equilibrium. Then either $u = Y^{\{i\}} \not\in
[0, U^{I_n\setminus\{k\}}]$ for some $i\in I_n\setminus\{k\}$ or,
by the definition of $U$, $u\in \gamma_j\cap[0,
U^{I_n\setminus\{k\}}]$ for some $j\in I_n\setminus\{k\}$. In the
former case, $U_i<a^{-1}_{ii}$. By $(\ref{3.1})$, $u = Y^{\{i\}}$
is below $\gamma_{\ell}$ for all $\ell\in I_n\setminus\{i\}$ so
$u$ is below $\gamma_k$. In the latter case, $u$ is below
$\gamma_k$ by $(C_k)$.}
\end{remark}

In some situation, $(\ref{3.1})$ has a simpler form. To see this,
suppose $Y^{\{i\}}$ is below $\gamma_j$ for all $j\in
I_n\setminus\{i\}$, i.e.
\begin{equation}
a_{ii} >a_{ji}, \quad j\in I_n\setminus\{i\}. \label{3.2}
\end{equation}
For each $j\in I_n\setminus\{i\}$, suppose also that either
\begin{equation}
a_{ik} \geq a_{jk}, \quad k\in I_n\setminus\{i\} \label{3.3}
\end{equation}
or $a_{jj}>a_{ij}$ and $\gamma_i\cap\gamma_j\cap(\cap_{\ell\in
I_n\setminus\{i, j\}}\pi_{\ell}) = \{E^{(ij)}\}$ which is below
$\gamma_{\ell}$ for all $\ell\in I_n\setminus\{i, j\}$. The latter
can be written
\begin{equation} \label{3.4}
a_{jj}>a_{ij}, \quad \frac{a_{\ell i}(a_{jj}-a_{ij})+a_{\ell
j}(a_{ii}-a_{ji})}{a_{ii}a_{jj} - a_{ij}a_{ji}} < 1, \quad \ell
\in I_n\setminus\{i, j\}.
\end{equation}
Then $(\ref{3.1})$ can be simplified as
\begin{equation} \label{3.5}
U_i = \max\left\{0,
\frac{a_{jj}-a_{ij}}{a_{ii}a_{jj}-a_{ij}a_{ji}}: j\in
I_n\setminus\{i\}, a_{jj}> a_{ij}\right\}.
\end{equation}
If condition $(C_k)$ holds for all $k\in I_n$, then, by remark
\ref{Re3.1}, $(\ref{3.2})$ and $(\ref{3.4})$ hold for all $i, j\in
I_n$ with $i\not=j$ so $U$ can be defined by
\begin{equation} \label{3.6}
U_i = \max\left\{\frac{a_{jj}-a_{ij}}{a_{ii}a_{jj} -
a_{ij}a_{ji}}: j\in I_n\setminus\{i\}\right\}, \quad i\in I_n.
\end{equation}

\begin{lemma} \label{Le3.2}
$(i)$ For any $J\subset I_n$ with $|J| >1$, condition $(C_i)$ for
all $i\in J$ implies
\begin{equation} \label{3.7}
\max\left\{0, \frac{a_{ij}}{a_{jj}}(1-\alpha_jU^{I_n\setminus\{i,
j\}})\right\} < 1-\alpha_iU^{I_n\setminus\{i, j\}}
\end{equation}
for all $i, j\in J$ with $i\not=j$. $(ii)$ Conversely, for a fixed
$i\in I_n$, $(\ref{3.7})$ for all $j\in I_n\setminus\{i\}$ implies
condition $(C_i)$. $(iii)$ Condition $(C_i)$ holds for all $i\in
I_n$ if and only if $(\ref{3.7})$ holds for all $i, j\in I_n$ with
$i\not=j$.
\end{lemma}

\textit{Proof.} $(i)$ Suppose $\alpha_iU^{I_n\setminus\{i, j\}}
\geq 1$ for some $i, j\in J$ with $i\not=j$. Then $[0,
U^{I_n\setminus\{i\}}]\cap\gamma_i\supset [0, U^{I_n\setminus\{i,
j\}}]\cap\gamma_i \not=\emptyset$. Since $[0, U^{I_n\setminus\{i,
j\}}]\cap\gamma_i$ is above $\gamma_j$ by $(C_i)$ whereas it is
below $\gamma_j$ by $(C_j)$, this is impossible. Therefore, we
must have
\begin{equation} \label{3.8}
\alpha_iU^{I_n\setminus\{i, j\}} <1, \quad i, j\in J, i\not=j.
\end{equation}
For any $i, j\in J$ with $i\not=j$, we assume $I_n\setminus\{i,
j\} \not=\emptyset$ as otherwise $(\ref{3.7})$ holds obviously. If
we can show that
\begin{equation} \label{3.9}
\alpha_jU^{I_n\setminus\{i\}} \geq 1
\end{equation}
for all $j\in J\setminus\{i\}$, then $(\ref{3.8})$ and
$(\ref{3.9})$ result in $\gamma_j\cap[U^{I_n\setminus\{i, j\}},
U^{I_n \setminus\{i\}}] = \{y\}$ with $y^{I_n\setminus\{j\}} =
U^{I_n\setminus\{i, j\}}$ and $y_j =
a^{-1}_{jj}(1-\alpha_jU^{I_n\setminus\{i, j\}})$. By condition
$(C_i)$ $y$ is below $\gamma_i$, so
\[
\alpha_iy = \alpha_iU^{I_n\setminus\{i, j\}} +
a_{ij}a^{-1}_{jj}(1-\alpha_jU^{I_n\setminus\{i, j\}}) <1
\]
and $(\ref{3.7})$ follows from this and $(\ref{3.8})$. For any
$j\in J\setminus\{i\}$, $(\ref{3.9})$ is obvious if $U_j =
a^{-1}_{jj}$. If $U_j < a^{-1}_{jj}$, then, for any $k\in
I_n\setminus\{i, j\}$, $Y^{\{j\}}$ is below $\gamma_k$ by the
definition of $U_j$. On the other hand, $Y^{\{k\}}$ is below
$\gamma_j$ by $(C_j)$ and remark \ref{Re3.1}. Consequently,
$\gamma_j\cap\gamma_k\cap(\cap_{\ell\in I_n\setminus\{j,
k\}}\pi_{\ell}) = \{E^{(jk)}\}$. Then $(\ref{3.9})$ follows from
$\alpha_jE^{(jk)} = 1$ and $E^{(jk)}\in [0,
U^{I_n\setminus\{i\}}]$ by the definition of $U$.

$(ii)$ If $\alpha_iU^{I_n\setminus\{i\}} <1$ then
$U^{I_n\setminus\{i\}}$ is below $\gamma_i$ so condition $(C_i)$
holds. Now suppose $\alpha_iU^{I_n\setminus\{i\}} \geq 1$. Then,
for each $j\in I_n\setminus\{i\}$, since
$\alpha_iU^{I_n\setminus\{i, j\}} <1$ by $(\ref{3.7})$, we have
$U_j>0$ and $Z^{(j)} \in [U^{I_n\setminus\{i, j\}},
U^{I_n\setminus\{i\}}]$ such that $\alpha_iZ^{(j)} = 1$. It can be
shown that
\[
\gamma_i\cap[0, U^{I_n\setminus\{i\}}] = \biggl\{\sum_{j\in
I_n\setminus\{i\}}d_jZ^{(j)}: d_j\geq 0, \sum_{j\in
I_n\setminus\{i\}}d_j =1\biggr\}.
\]
Hence, to show that condition $(C_i)$ is met, we need only verify
\begin{equation}\label{3.10}
\alpha_mZ^{(j)} >1
\end{equation}
for each fixed $j\in I_n\setminus\{i\}$ and all $m\in
I_n\setminus\{i\}$. Note that $Z^{(j)}_j =
a^{-1}_{ij}(1-\alpha_iU^{I_n\setminus\{i, j\}})$. It then follows
from this and $(\ref{3.7})$ that
\begin{equation} \label{3.11}
\alpha_jZ^{(j)} = \alpha_jU^{I_n\setminus\{i, j\}} +
a_{jj}a^{-1}_{ij}(1 - \alpha_iU^{I_n\setminus\{i, j\}}) >1
\end{equation}
for all $j\in I_n\setminus\{i\}$. For any $m\in I_n\setminus\{i,
j\}$, if $a_{mj}= 0$ then $\alpha_mZ^{(j)} =
\alpha_mU^{I_n\setminus\{i\}} \geq \alpha_mZ^{(m)} >1$ by
$(\ref{3.11})$. If $a_{mj}>0$ and $U_m = a^{-1}_{mm}$, we also
have $\alpha_mZ^{(j)} \geq a_{mm}U_m +a_{mj}Z^{(j)}_j >a_{mm}U_m
=1$. Since $a_{mj} =0$ implies $U_m = a^{-1}_{mm}$, we have
actually verified $(\ref{3.10})$ when $U_m = a^{-1}_{mm}$. So our
remaining task is to show $(\ref{3.10})$ under the assumption $U_m
< a^{-1}_{mm}$, which implies $a_{mm}>a_{km}$ and $a_{mk}>0$ for
all $k\in I_n\setminus\{m\}$ by the definition of $U_m$.

(A) Suppose $I_n\not = \{i, j, m\}$. Then, for $k\in
I_n\setminus\{i, j, m\}$, we have either $0<a_{kk}\leq a_{mk}$ or
$0< a_{mk}< a_{kk}$. In the former case, $U_k = a^{-1}_{kk}$ so
$\alpha_mZ^{(j)} > a_{mk}U_k \geq 1$. In the latter case,
$\gamma_m\cap\gamma_k\cap(\cap_{\ell\in I_n\setminus\{m,
k\}}\pi_{\ell}) = \{E^{(mk)}\}$. Since $E^{(mk)}\leq U^{\{m, k\}}$
and $\alpha_mE^{(mk)} =1$, we have $\alpha_mZ^{(j)} >
\alpha_mU^{\{m, k\}} \geq \alpha_mE^{(mk)} = 1$.

(B) Suppose $I_n = \{i, j, m\}$. Since
$\alpha_jU^{I_n\setminus\{i\}} = a_{jj}U_j + a_{jm}U_m \geq
\alpha_jZ^{(j)} >1$ by $(\ref{3.11})$ and $a_{jj}U_j\leq 1$, we
must have $a_{jm}U_m >0$ so $a_{mm}>a_{jm}>0$. If $a_{mj}\geq
a_{jj}$ then $\alpha_mZ^{(j)} \geq \alpha_jZ^{(j)} >1$. Otherwise,
we have $a_{jj} >a_{mj} >0$ so $\gamma_j\cap\gamma_m\cap\pi_i =
\{E^{(jm)}\}$ with
\[
(E^{(jm)}_j, E^{(jm)}_m) = \biggl(\frac{a_{mm} -
a_{jm}}{a_{jj}a_{mm}-a_{jm}a_{mj}}, \frac{a_{jj} -
a_{mj}}{a_{jj}a_{mm}-a_{jm}a_{mj}}\biggr).
\]
By the definition of $U$ we have $E^{(jm)} \leq U^{\{j, m\}}$.
Since $E^{(jm)}_m \leq U_m$ is equivalent to
$a^{-1}_{jj}(1-a_{jm}U_m) \geq a^{-1}_{mj}(1-a_{mm}U_m)$ and
$(\ref{3.11})$ gives $a^{-1}_{ij}(1-a_{im}U_m) > a^{-1}_{jj}(1-
a_{jm}U_m)$, from these we obtain
\begin{eqnarray*}
\alpha_mZ^{(j)} &=& a_{mm}U_m + a_{mj}a^{-1}_{ij}(1-a_{im}U_m) \\
    &>& a_{mm}U_m + a_{mj}a^{-1}_{mj}(1- a_{mm}U_m) \\
    &=& 1.
\end{eqnarray*}
Therefore, $(\ref{3.10})$ holds for all $m\in I_n\setminus\{i\}$.

$(iii)$ This follows from $(i)$ and $(ii)$.

\section{Restatement of the main results}

With the preparation in \S3 we are now in a position to restate
the theorems and corollaries given in \S2.

\begin{theorem} \label{Th4.1}
Assume that $(\ref{1.1})$ satisfies $(\ref{3.2})$, $(\ref{3.4})$
and $(\ref{3.7})$ for all $i, j\in I_n$ with $i\not= j$, where $U$
is given by $(\ref{3.6})$. Then $(\ref{1.1})$ has an equilibrium
$x^*\in \textup{int}{\R}^n_+$ with $x^*\leq U$ that is a global
attractor.
\end{theorem}

\begin{corollary} \label{Co4.2}
Under the conditions of theorem \ref{Th4.1}, for each nonempty set
$J\subset I_n$, $(\ref{1.1})$ has an equilibrium $x_J\in [0, U^J]$
with $(x_J)_i >0$ for $i\in J$ that is a global attractor within
the interior of $\cap_{\ell\in I_n\setminus J}\pi_{\ell}$ viewed
as ${\R}^{|J|}_+$. Moreover, $(\ref{1.1})$ has $2^n$ equilibria in
${\R}^n_+$.
\end{corollary}

From lemma \ref{Le3.2} we see that the conditions of theorem
\ref{Th4.3} and corollary \ref{Co4.4} given below are slightly
stronger than those of theorem \ref{Th2.5} and corollary
\ref{Co2.8} respectively.

\begin{theorem} \label{Th4.3}
Let $J\subset I_n$ with $J\not=\emptyset$. Assume that
$(\ref{3.7})$ holds for each $i\in J$ and every $j\in
I_n\setminus\{i\}$, where each component $U_i$ of $U$ is defined
either by $(\ref{3.1})$ or, if $(\ref{3.2})$ and $(\ref{3.4})$
holds for all $j\in I_n\setminus\{i\}$, by $(\ref{3.6})$. Then
$(\ref{1.1})$ has an equilibrium $x^*\in [0, U^J]$ with $x^*_j >0$
for $j\in J$. In addition, if $\alpha_jx^* = 1$ for all $j\in
I_n\setminus J$, then $x^*$ is a global attractor. Moreover,
``$\alpha_jx^* = 1$'' can be relaxed to ``$\alpha_jx^* \geq 1$''
when $U$ is replaced by $Y$ given by $(\ref{2.1})$.
\end{theorem}

\begin{corollary} \label{Co4.4}
Assume that $(\ref{3.7})$ and $a_{ii} = a_{ji}$ hold for some
$i\in I_n$ and all $j\in I_n\setminus\{i\}$, where $U$ is defined
by $(\ref{3.1})$. Then $Y^{\{i\}}$ is a global attractor of
$(\ref{1.1})$. Moreover, ``$a_{ii} = a_{ji}$'' can be replaced by
``$a_{ii}\leq a_{ji}$'' if $U$ is replaced by $Y$.
\end{corollary}

\begin{theorem} \label{Th4.5} Let $J = \{j_1, j_2, \ldots,
j_m\}\subset I_n$ with $1\leq m <n$ and let $U\in {\R}^n_+$ be
defined by $U_i = 0$ for $i\in J$ and by $(\ref{3.1})$, with $I_n$
replaced by $I_n\setminus J$, for $i\in I_n\setminus J$. Assume
that
\begin{itemize}
\item[$(i)$] for each $\ell\in I_m$ and every $j, k\in
I_n\setminus\{j_1, j_2, \ldots, j_{\ell}\}$,
\begin{equation} \label{4.1}
a_{kj_{\ell}} < a_{j_{\ell}j_{\ell}}, \quad a_{kj} \leq
a_{j_{\ell}j};
\end{equation}
\item[$(ii)$] for all $i, j\in I_n\setminus J$ with $i\not= j$ (if
$n> m+1$),
\begin{equation} \label{4.2}
\max\biggl\{0, \frac{a_{ij}}{a_{jj}}(1 - \alpha_jU^{(I_n\setminus
J)\setminus\{i, j\}})\biggr\} < 1 - \alpha_iU^{(I_n\setminus
J)\setminus\{i, j\}}.
\end{equation}
\end{itemize}
Then $(\ref{1.1})$ has an equilibrium $x^*\in [0, U^{I_n\setminus
J}]$ with $x^*_j >0$ for $j\in I_n\setminus J$ that is a global
attractor.
\end{theorem}

\begin{corollary}\label{Co4.6} Assume that $(\ref{1.1})$
satisfies $(\ref{2.5})$. Then $Y^{\{n\}}$ is a global attractor.
\end{corollary}

\begin{remark} \label{Re4.7}
\textup{By remark \ref{Re3.1} and theorem \ref{Th2.1}, we are
tempted to make the following conjecture as an improvement of
theorem \ref{Th2.1}:
\begin{quotation}
\noindent System $(\ref{1.1})$ has an equilibrium $x^*\in
\textup{int}{\R}^n_+$ that is a global attractor if, for each
$i\in I_n$, every equilibrium of $(\ref{1.1})$ on $\pi_i$ is below
$\gamma_i$.
\end{quotation}
This is obviously true when $n=2$. When $n= 3$, the condition
ensures that each $Y^{\{i\}}$ is a local repellor and, by lemma
\ref{Le5.1} given in the next section,
$\gamma_1\cap\gamma_2\cap\gamma_3 = \{x^*\} \subset
\textup{int}{\R}^3_+$. Then, by a result given in \cite{DrZe},
$x^*$ is a global attractor. Hence, it is fairly reasonable to
make the above conjecture for general $n$-dimensional system
$(\ref{1.1})$. However, further investigation is needed to clarify
it.}
\end{remark}

\section{The proofs of theorems \ref{Th2.5} and
\ref{Th2.10}}

The method of ultimate contracting cells (see \cite[Theorem
2]{Ho2} or \cite[Lemma 1]{Ho3}) will be used in the proofs of the
theorems though it is not explicitly stated here. Our first lemma
establishes the existence of a unique equilibrium $x^*\in
\textup{int}{\R}^n_+$ under conditions weaker than $(C_k)$ for all
$k\in I_n$.

\begin{lemma}\label{Le5.1}
Assume that, for each $i\in I_n$, every equilibrium of
$(\ref{1.1})$ on $\pi_i$ is below $\gamma_i$. Then there is a
unique equilibrium $x^*\in \textup{int}{\R}^n_+$ with $x^*\leq Y$.
\end{lemma}

\begin{figure}
\begin{picture}(300, 120)(80, 10) \thicklines \put(110,
30){\line(1, 0){160}} \put(110, 30){\line(1, 1){80}} \put(215,
30){\line(-3, 4){45}} \put(236, 30){\line(-2, 1){84}} \put(100,
20){$v^*$} \put(205, 20){$w^*$} \put(225, 20){$u(\delta_1)$}
\put(272, 25){$\pi_{\ell}$} \put(142, 77){$u^*$} \put(192,
114){$\pi_{k+1}$} \put(180, 85){$\gamma_{k+1}$} \put(225,
40){$\gamma_{\ell}$} \put(190, 40){$x^*$} \put(160, 10){Figure
5.1} \put(250, 100){$S$}
\end{picture}
\end{figure}

\textit{Proof.} Denoting the statement of lemma \ref{Le5.1} by
$P(n)$, we show the truth of $P(n)$ by induction. The truth of
$P(n)$ is trivial for $n=1$ and $n=2$. Supposing that $P(n)$ is
true for some $k\geq 2$ and all $n\in I_k$, we now show the truth
of $P(k+1)$.

By the assumption, for each $i\in I_{k+1}$, every equilibrium of
$(\ref{1.1})$ on $\pi_i$ is below $\gamma_i$. In particular, for
each $i\in I_k$, every equilibrium of $(\ref{1.1})$ on
$\pi_i\cap\pi_{k+1}$ is below $\gamma_i$ and $\gamma_{k+1}$.
Viewing $(\ref{1.1})$ on $\pi_{k+1}$ as a $k$-dimensional system
on ${\R}^k_+$ and using the truth of $P(k)$, we have a unique
equilibrium $u^*$ of $(\ref{1.1})$ on $\pi_{k+1}$ satisfying $0<
u^*_i\leq a^{-1}_{ii}$ for $i\in I_k$. This implies that the
solution of the system
\[
\alpha_iu =1, \quad i\in I_k
\]
can be written $u(\delta)$ with $u(0) = u^*$, $u_{k+1}(\delta) =
\delta$ and $u(\delta)-u^*$ is linear in $\delta$. Thus,
$u(\delta) \in \textup{int}{\R}^{k+1}_+$ holds for sufficiently
small $\delta >0$.

$(i)$ If $u(a^{-1}_{k+1 k+1})\in \textup{int}{\R}^{k+1}_+$ then
$u(\delta)\in \textup{int}{\R}^{k+1}_+$ for all $\delta \in (0,
a^{-1}_{k+1 k+1}]$ and
\[
\alpha_{k+1}u(a^{-1}_{k+1 k+1}) \geq a_{k+1 k+1}a^{-1}_{k+1 k+1} =
1.
\]
As $\alpha_{k+1}u(0) = \alpha_{k+1}u^* <1$ by the assumption,
there is a $\delta_0\in (0, a^{-1}_{k+1 k+1}]$ such that
$\alpha_{k+1}u(\delta_0) =1$. Then $u(\delta_0) \in
\textup{int}{\R}^{k+1}_+$ is the required equilibrium.

$(ii)$ If $u(a^{-1}_{k+1 k+1})\not\in \textup{int}{\R}^{k+1}_+$,
then there is a $\delta_1\in (0, a^{-1}_{k+1 k+1}]$ such that
$u(\delta) \in \textup{int}{\R}^{k+1}_+$ for all $\delta\in (0,
\delta_1)$ but $u_{\ell}(\delta_1)=0$ for some $\ell \in I_k$.
Thus $\cap_{i\in I_k}\gamma_i$ is a line segment joining $u^*$
with $u(\delta_1)$. Again, by the assumption and the truth of
$P(k-1)$ and $P(k)$, $(\ref{1.1})$ has a unique equilibrium $v^*$
on $\pi_{\ell}\cap\pi_{k+1}$ satisfying $0< v^*_i\leq a^{-1}_{ii}$
for $i\in I_k\setminus\{\ell\}$ and a unique equilibrium $w^*$ on
$\pi_{\ell}$ with $0< w^*_j\leq a^{-1}_{jj}$ for $j\in
I_{k+1}\setminus\{\ell\}$. Consequently, the set $S = \cap_{i\in
I_k\setminus\{\ell\}}\gamma_i$ is a two-dimensional plane and
$\cap_{i \in I_{k+1}\setminus\{\ell\}}\gamma_i$ is a line segment
starting from $w^*$ on $S$ (see Figure 5.1). Note that
$S\cap\pi_{\ell}$ is a line segment containing $v^*$, $w^*$ and
$u(\delta_1)$, and $S\cap\pi_{k+1}$ is a line segment containing
$v^*$ and $u^*$. Since $v^*$ is below both $\gamma_{\ell}$ and
$\gamma_{k+1}$, $u^*$ is below $\gamma_{k+1}$ and $w^*$ is below
$\gamma_{\ell}$, we must have $\gamma_{\ell}\cap\gamma_{k+1}\cap S
=\{x^*\} \subset \textup{int}{\R}^n_+$ with $x^*$ the required
equilibrium.

Therefore, we have shown the truth of $P(k+1)$ when $P(n)$ is true
for all $n\in I_k$. By induction, $P(n)$ is true for all $n\geq
1$.

The next lemma gives an ultimate upper bound of any solution of
$(\ref{1.1})$ in terms of a given ultimate lower bound (see
\cite[Lemma 3.2]{Ho1}).

\begin{lemma}\label{Le5.2}
Let $\hat{x}\in{\R}^n_+$ satisfy $\hat{x}_i = 0$ if
$\alpha_i\hat{x}>1$ for any $i\in I_n$ and let
$\hat{y}\in{\R}^n_+$ be given by
\begin{equation} \label{5.1}
\hat{y}_j = \max\{0, \hat{x}_j + a^{-1}_{jj}(1-\alpha_j\hat{x})\},
\quad j\in I_n.
\end{equation}
If a solution of $(\ref{1.1})$ in $\textup{int}{\R}^n_+$ satisfies
$\underline{\lim}_{t\to\infty}x(t) \geq \hat{x}$, then it also
satisfies $\overline{\lim}_{t\to\infty}x(t) \leq \hat{y}$.
\end{lemma}

However, we shall need a refined ultimate upper bound based on
$(\ref{5.1})$. For the above $\hat{x}$ and some $j\in I_n$, if we
know that $x_j(t)\to\hat{x}_j$ as $t\to\infty$, then we take
$\hat{y}_j = \hat{x}_j$ rather than $(\ref{5.1})$. Now for
$\hat{y}\geq \hat{x}$ with $\hat{y}_j = \hat{x}_j$ or $\hat{y}_j$
given by $(\ref{5.1})$ for each $j\in I_n$, let $J_0= \{i\in I_n:
\hat{x}_i <\hat{y}_i\}$ and define $\hat{z}\in [\hat{x}, \hat{y}]$
by $(\ref{5.2})$ for each $i\in I_n$,
\begin{equation} \label{5.2}
\hat{z}_i = \left\{\begin{array}{l} \hat{y}_i \quad \textrm{if}\;
J_0= \{i\} \;\textrm{or}\; i\not\in J_0 \; \textrm{or, for some}\;
j\in J_0\setminus\{i\}, \\
\qquad a_{ij}=0 \; \textrm{or}\;
\alpha_j\hat{x}^{I_n\setminus\{i\}} + a_{ji}\hat{y}_i \geq 1; \\
\hat{x}_i \quad \textrm{if} \; i\in J_0\not=\{i\},
\alpha_j\hat{x}^{I_n\setminus\{i\}} +a_{ji}\hat{y}_i <1 \;
\textrm{and}\; \alpha_jy\leq 1 \\
\qquad \textrm{for all}\; y\in\gamma_i\cap[\hat{x}, \hat{y}] \;
\textrm{and}\; j\in J_0\setminus\{i\}; \\
\max\{u^{(ijk)}_i: \gamma_i\cap\gamma_k\cap[\hat{x}, H^{(ij)}] =
\{u^{(ijk)}\}, j, k\in J_0\setminus\{i\}\} \\
\qquad \textrm{otherwise},
\end{array}\right.
\end{equation}
where $H^{(ij)} = \hat{x}^{I_n\setminus\{i, j\}} + \hat{y}^{\{i,
j\}}$.

\begin{lemma}\label{Le5.3}
Assume that a solution $x(t)$ of $(\ref{1.1})$ in
\textup{int}${\R}^n_+$ satisfies
\[
\hat{x} \leq \underline{\lim}_{t\to\infty}x(t) \leq
\overline{\lim}_{t\to\infty}x(t) \leq \hat{y}
\]
for some $\hat{x}, \hat{y} \in \textup{int}{\R}^n_+$, where, for
any $i, j \in I_n$, $\hat{x}_i = 0$ if $\alpha_i\hat{x} > 1$ and
$\hat{y}_j$ is given by $(\ref{5.1})$ if $\hat{y}_j > \hat{x}_j$.
Then, for $\hat{z}\geq\hat{y}$ defined by $(\ref{5.2})$, $x(t)$
also satisfies $\overline{\lim}_{t\to\infty}x(t) \leq \hat{z}$.
\end{lemma}

\textit{Proof.} We only need to show
$\overline{\lim}_{t\to\infty}x_i(t) \leq \hat{z}_i$ if $\hat{z}_i
< \hat{y}_i$ for a fixed $i\in I_n$. This can be achieved by
showing that, for every $r_0\in (\hat{z}_i, \hat{y}_i]$, there is
an $r_1\in [\hat{z}_i, r_0)$ such that
\begin{equation} \label{5.3}
\overline{\lim}_{t\to\infty}x_i(t) \leq r_0
\end{equation}
implies
\begin{equation}\label{5.4}
\overline{\lim}_{t\to\infty}x_i(t) \leq r_1.
\end{equation}
From $\hat{z}_i<\hat{y}_i$ and $(\ref{5.2})$ we know that $i\in
J_0\not=\{i\}$, $a_{ij}>0$ and
$\alpha_j\hat{x}^{I_n\setminus\{i\}} + a_{ji}\hat{y}_i <1$, i.e.
$\hat{x}^{I_n\setminus\{i\}} + \hat{y}^{\{i\}}$ is below
$\gamma_j$, for all $j\in J_0\setminus\{i\}$. By $(\ref{5.1})$,
$\hat{x}^{I_n\setminus\{k\}} + \hat{y}^{\{k\}}$ is on $\gamma_k$
for all $k\in J_0$. Thus, for each $j\in J_0\setminus\{i\}$,
$H^{(ij)} = \hat{x}^{I_n\setminus\{i, j\}}+ \hat{y}^{\{i, j\}}$ is
on or above $\gamma_j$ and $\hat{x}^{I_n\setminus\{i\}}
+\hat{y}^{\{i\}}$ is below $\gamma_k$ for every $k\in
J_0\setminus\{i\}$. Now suppose $(\ref{5.3})$ holds. The line
segment $L_{r_0} = \{x\in[\hat{x}, H^{(ij)}]: x_i = r_0\}$ must
intersect a plane $\gamma_k$ for some $k\in J_0\setminus\{i\}$ at
$Q^{(j)}$ (see Figure 5.2) such that $Q^{(j)}_i = r_0$,
$Q^{(j)}_{\ell} =\hat{x}_{\ell}$ for $\ell\in I_n\setminus\{i,
j\}$ and
\[
Q^{(j)}_j = \min\{a^{-1}_{\ell j}(1 -
\alpha_{\ell}\hat{x}^{I_n\setminus\{i, j\}} - a_{\ell i}r_0):
\ell\in J_0\setminus\{i\}, a_{\ell j}\not= 0\}.
\]
By the definition of $\hat{z}_i$ and $r_0>\hat{z}_i$, we see that
$Q^{(j)}_j >\hat{x}_j$ and $Q^{(j)}$ is above $\gamma_i$ but on or
below $\gamma_{\ell}$ for every $\ell\in J_0\setminus\{i\}$. Put
$\gamma(\varepsilon) = \{x\in{\R}^n_+: \beta x = \varepsilon\}$,
where $\beta = (\beta_1, \beta_2, \ldots, \beta_n)$ with $\beta_j
=0$ for $j\in \{i\}\cup(I_n\setminus J_0)$ and
\[
\beta_j = \frac{(Q^{(j)}_j - \hat{x}_j)^{-1}}{1 + \sum_{\ell\in
J_0\setminus\{i\}}(Q^{(\ell)}_{\ell}
-\hat{x}_{\ell})^{-1}\hat{x}_{\ell}}
\]
for $j\in J_0\setminus\{i\}$. Then the $x_i$-axis is below and
parallel to the plane $\gamma(\varepsilon)$ for each $\varepsilon
>0$, $[\hat{x}, \hat{x}^{I_n\setminus\{i\}}+\hat{y}^{\{i\}}]$ is
below and parallel to $\gamma(1)$, and $Q^{(j)}$ is on $\gamma(1)$
for $j\in J_0\setminus\{i\}$. Moreover, for any
$\varepsilon\in(\beta\hat{x}, 1)$, the set
\begin{equation}\label{5.5}
A(\varepsilon) = \{x\in[\hat{x}, \hat{y}]: x_i\leq r_0, \beta
x\leq \varepsilon\}
\end{equation}
is below $\gamma_j$ for every $j\in J_0\setminus\{i\}$. We now
show the inequality
\begin{equation}\label{5.6}
\underline{\lim}_{t\to\infty}\beta x(t) \geq 1
\end{equation}
by assuming its falsity, i.e. $\underline{\lim}_{t\to\infty}\beta
x(t) = \delta$ for some $\delta<1$. Then there is a sequence
$\{t_m\}$ satisfying $t_m\to\infty$ as $m\to\infty$ and $\beta
x(t_m) < (1+\delta)/2$ for all $m\geq 1$. Since
$\lim_{t\to\infty}x_{\ell}(t) = \hat{x}_{\ell}$ for $\ell\in
I_n\setminus J_0$, $(\ref{5.3})$ holds and $A(\frac{1+\delta}{2})$
given by $(\ref{5.5})$ is a closed set, we can assume the
existence of $T>0$, $[\tilde{x}, \tilde{y}] \supseteq [\hat{x},
\hat{y}]$ and $\tilde{r}_0 >r_0$ such that $x(t)\in [\tilde{x},
\tilde{y}]$ and $x_i(t) \leq \tilde{r}_0$ for all $t\geq T$ and
the set
\[
\tilde{A}((1+\delta)/2) = \{x\in [\tilde{x}, \tilde{y}]: x_i\leq
\tilde{r}_0, \beta x \leq (1+\delta)/2\}
\]
is below $\gamma_j$ for all $j\in J_0\setminus\{i\}$. Thus, for
any $t\geq T$, $x(t)\in \tilde{A}(\frac{1+\delta}{2})$ implies
$\beta x'(t) >0$. Without loss of generality, we may also assume
that $x(t_m)\in\tilde{A}(\frac{1+\delta}{2})$, so that $\beta
x'(t_m)
>0$, for all $m\geq 1$. Then, for each $m\geq 1$, $\beta x(t)$ is
increasing for $t\geq t_m$ as long as $x(t)\in
\tilde{A}(\frac{1+\delta}{2})$. This, together with $\beta
x(t_m)<(1+\delta)/2$ for all $m\geq 1$, results in $x(t)\in
\tilde{A}(\frac{1+\delta}{2})$ for all $t\geq t_1$. By the
closeness of $\tilde{A}(\frac{1+\delta}{2})$ there is an
$\varepsilon\in (0, 1)$ such that $\alpha_jx\leq\varepsilon$ for
all $x\in\tilde{A}(\frac{1+\delta}{2})$ and $j\in
J_0\setminus\{i\}$. Hence, from $(\ref{1.1})$,
\[
x_j(t) \geq x_j(t_1) e^{b_j(1-\varepsilon)(t-t_1)}
\]
for $t\geq t_1$ and $j\in J_0\setminus\{i\}$. This leads to $\beta
x(t) \to\infty$ as $t\to\infty$, a contradiction to
$\underline{\lim}_{t\to\infty}\beta x(t) =\delta$. Therefore, we
must have $(\ref{5.6})$.

Note that, for each $j\in J_0\setminus\{i\}$,
$\gamma(1)\cap[\hat{x}, H^{(ij)}]$ is a line segment parallel to
$[\hat{x}, \hat{x}^{I_n\setminus\{i\}}+\hat{y}^{\{i\}}]$ and
$Q^{(j)}$ is on it. So we let $\gamma_i\cap\gamma(1)\cap[\hat{x},
H^{(ij)}] = \{P^{(j)}\}$ or, if it is empty, let
$\gamma(1)\cap[\hat{x}, \hat{x}^{I_n\setminus\{j\}} +
\hat{y}^{\{j\}}] = \{P^{(j)}\}$, and let
\[
r_1 = \max\{P^{(j)}_i: j\in J_0\setminus\{i\}\}.
\]
Then $r_1\in [\hat{z}_i, r_0)$ as each $Q^{(j)}$ is above
$\gamma_i$. We observe that the closed set $S(\eta) =
\{x\in[\hat{x}, \hat{y}]: x_i\geq \eta, \beta x \geq 1\}$ is above
$\gamma_i$ for every $\eta >r_1$. Then, from this observation and
$(\ref{5.6})$ we see that $x(t)$ satisfies $x_i(t) < \eta$ for
each $\eta> r_1$ and all sufficiently large $t$. Thus
$(\ref{5.4})$ holds.

\begin{figure}
\begin{picture}(350, 180)(100, 10) \thicklines \put(120,
30){\line(1, 0){200}} \put(120, 30){\line(0, 1){140}} \put(120,
170){\line(1, 0){200}} \put(320, 30){\line(0, 1){140}} \put(320,
30){\line(-1, 1){140}} \put(120, 90){\line(1, 0){200}} \put(120,
150){\line(3, -1){200}} \put(300, 30){\line(0, 1){140}} \put(240,
30){\line(0, 1){140}} \put(120, 170){\line(5, -2){200}} \put(110,
20){$\hat{x}$} \put(220,
20){$\hat{x}^{I_n\setminus\{i\}}+\hat{z}^{\{i\}}$} \put(305,
20){$\hat{x}^{I_n\setminus\{i\}}+\hat{y}^{\{i\}}$} \put(100,
172){$\hat{x}^{I_n\setminus\{j\}}+\hat{y}^{\{j\}}$} \put(210,
172){$\hat{x}^{I_n\setminus\{i, j\}}+\hat{z}^{\{i, j\}}$}
\put(320, 172){$H^{(ij)}$} \put(140, 92){$\gamma(1)\cap[\hat{x},
H^{(ij)}]$} \put(130, 147){$\gamma_k$} \put(150, 159){$\gamma_j$}
\put(301, 130){$L_{r_0}$} \put(250, 80){$P^{(j)}$} \put(301,
80){$Q^{(j)}$} \put(200, 151){$\gamma_i$} \put(200, 5){Figure 5.2}
\end{picture}
\end{figure}

From lemmas \ref{Le5.2} and \ref{Le5.3} we obtain an ultimate
upper bound for all solutions of $(\ref{1.1})$ in int${\R}^n_+$.

\begin{lemma}\label{Le5.4}
Every solution of $(\ref{1.1})$ in \textup{int}${\R}^n_+$
satisfies
\[
\overline{\lim}_{t\to\infty}x(t) \leq U
\]
for $U$ given by $(\ref{3.1})$. If $U_{j_1} = 0$ for some $j_1\in
I_n$, then there exists $J = \{j_1, j_2, \ldots, j_m\}\subset I_n$
with $1\leq m < n$ such that $(\ref{4.1})$ holds for each $\ell\in
I_m$ and every $j, k\in I_n\setminus\{j_1, j_2, \ldots,
j_{\ell}\}$. Moreover, every solution of $(\ref{1.1})$ in
\textup{int}${\R}^n_+$ satisfies
\[
\overline{\lim}_{t\to\infty}x(t) \leq V,
\]
where $V_k = 0$ for $k\in J$ and $V_i\in (0, a^{-1}_{ii}]$ is
given by $(\ref{3.1})$ with the replacement of $I_n$ by
$I_n\setminus J$ for $i\in I_n\setminus J$.
\end{lemma}

\textit{Proof.} By lemma~\ref{Le5.2}, every solution of
$(\ref{1.1})$ in \textrm{int}${\R}^n_+$ satisfies
\[
0 = \hat{x} \leq
\underline{\lim}_{t\to\infty}x(t) \leq
\overline{\lim}_{t\to\infty}x(t) \leq \hat{y} = Y.
\]
Note that $\hat{z}$ given by $(\ref{5.2})$ with $[\hat{x},
\hat{y}] = [0, Y]$ coincides with $U$ defined by $(\ref{3.1})$.
Then, by lemma~\ref{Le5.3}, every solution of $(\ref{1.1})$ in
\textrm{int}${\R}^n_+$ satisfies
\begin{equation}\label{5.7}
0 = \hat{x} \leq \underline{\lim}_{t\to\infty}x(t) \leq
\overline{\lim}_{t\to\infty}x(t) \leq \hat{z}
\end{equation}
with $\hat{z} = U$. If $U_{j_1}=0$ for some $j_1\in I_n$, then,
with $J = \{j_1\}$, $(\ref{4.1})$ follows from $(\ref{3.1})$ for
$\ell = 1$ and all $j, k\in I_n\setminus\{j_1\}$. From the
observation $(ii)$ made before theorem \ref{Th2.10} we know that
$U = Y^{I_n\setminus\{J_1\}}$. Applying lemma~\ref{Le5.3} again
with $[\hat{x}, \hat{y}] = [0, Y^{I_n\setminus\{j_1\}}]$, we see
that every solution of $(\ref{1.1})$ in \textrm{int}${\R}^n_+$
satisfies $(\ref{5.7})$ with $\hat{z}$ given by $(\ref{5.2})$.
Note that, for each $j\in I_n\setminus\{j_1\}$, $\hat{z}_j$
coincides with $U_j$ given by $(\ref{3.1})$ after the replacement
of $I_n$ by $I_n\setminus\{j_1\}$. If $\hat{z}_j >0$ for all
$j\not = j_1$, then the conclusion follows with $J = \{j_1\}$ and
$V = \hat{z}$. Otherwise, if $\hat{z}_{j_2} = U_{j_2}=0$ for some
$j_2\in I_n\setminus\{j_1\}$, the conclusion follows from
repeating the above process a number of times.

The next lemma gives a refined ultimate lower bound based on the
given ultimate lower and upper bounds. In particular, it confirms
that condition $(C_k)$ guarantees the survival of the $k$th
species.

\begin{lemma}\label{Le5.5}
Assume that every solution of
$(\ref{1.1})$ in $\textup{int}{\R}^n_+$ satisfies
\begin{equation}\label{5.8}
\hat{x} \leq \underline{\lim}_{t\to\infty}x(t) \leq
\overline{\lim}_{t\to\infty}x(t) \leq \hat{u}
\end{equation}
for some $\hat{x}, \hat{u} \in{\R}^n_+$. Assume also that, for
some $i\in J_1 = \{j\in I_n: \hat{u}_j > \hat{x}_j\} \not=
\emptyset$, either $(A)$ $\hat{x}^{\{i\}} +
\hat{u}^{I_n\setminus\{i\}}$ is below $\gamma_i$ or $(B)$
$\hat{x}$ is below $\gamma_i$ but $\gamma_i\cap[\hat{x},
\hat{x}^{\{i\}} +\hat{u}^{I_n\setminus\{i\}}]$ is above $\gamma_j$
for every $j\in J_1\setminus\{i\}$. Then there is a $\delta_0>0$
such that
\begin{equation}\label{5.9}
\underline{\lim}_{t\to\infty}x_i(t) \geq \hat{x}_i +\delta_0.
\end{equation}
\end{lemma}

\textit{Proof.} Under $(A)$ and $(\ref{5.8})$ we show that
\[
\delta_0 = a^{-1}_{ii}(1- a_{ii}\hat{x}_i -
\alpha_i\hat{u}^{I_n\setminus\{i\}})
\]
meets the requirement of $(\ref{5.9})$. In fact, $\delta_0>0$
since $\hat{x}^{\{i\}}+\hat{u}^{I_n\setminus\{i\}}$ is below
$\gamma_i$. Moreover, for arbitrary $\delta\in (0, \delta_0)$, we
have $\alpha_i\hat{u}^{I_n\setminus\{i\}} + a_{ii}(\delta_0
+\hat{x}_i) = 1$ and $\alpha_i\hat{u}^{I_n\setminus\{i\}} +
a_{ii}(\delta +\hat{x}_i) < 1$. For each solution of $(\ref{1.1})$
in $\textup{int}{\R}^n_+$, since
$\overline{\lim}_{t\to\infty}x(t)^{I_n\setminus\{i\}} \leq
\hat{u}^{I_n\setminus\{i\}}$ by $(\ref{5.8})$, there is a $t_0>0$
and $\varepsilon\in (0, 1)$ such that
$\alpha_ix(t)^{I_n\setminus\{i\}} + a_{ii}(\delta +\hat{x}_i) \leq
\varepsilon < 1$ for all $t\geq t_0$. If $x_i(t) \leq \hat{x}_i +
\delta$ holds for all $t\geq t_0$ then
\[
\alpha_ix(t) \leq \alpha_ix(t)^{I_n\setminus\{i\}} + a_{ii}(\delta
+\hat{x}_i) \leq \varepsilon < 1
\]
for $t\geq t_0$, so $x_i(t) \geq
x_i(t_0)e^{b_i(1-\varepsilon)(t-t_0)} \to\infty$ as $t\to\infty$.
This is impossible. Since $x_i(\bar{t})\leq \hat{x}_i + \delta$
for any $\bar{t} \geq t_0$ implies the increase of $x_i$ for $t$
at the vicinity of $\bar{t}$, there must be a $t_1>t_0$ such that
$x_i(t) > \hat{x}_i+\delta$ for all $t\geq t_1$. Thus,
$\underline{\lim}_{t\to\infty}x_i(t) \geq \hat{x}_i +\delta$. Then
$(\ref{5.9})$ follows from this and the arbitrariness of
$\delta\in (0, \delta_0)$.

Under $(\ref{5.8})$ and $(B)$, the proof of \cite[Lemma 3.3]{Ho1}
can be copied to here almost verbatim after the replacement of
$[\hat{x}, \hat{y}]$ by $[\hat{x}, \hat{u}]$.

With the above preparation we are now able to embark on the proofs
of the theorems.

\vskip 5 pt \textit{Proof of theorem \ref{Th2.5}.} For $k\in J$,
since condition $(C_k)$ holds, by remark \ref{Re3.1} every
equilibrium of $(\ref{1.1})$ on $\pi_k\cap(\cap_{j\in I_n\setminus
J}\pi_j)$ is below $\gamma_k$. Restricting $(\ref{1.1})$ to
$\cap_{j\in I_n\setminus J}\pi_j$ and applying lemma \ref{Le5.1}
to it, we obtain a unique equilibrium $x^*$ on $\cap_{j\in
I_n\setminus J}\pi_j$ satisfying $0<x^*_k\leq a^{-1}_{kk}$ for
$k\in J$. By lemma \ref{Le5.4}, every solution of $(\ref{1.1})$ in
$\textup{int}{\R}^n_+$ satisfies
\begin{equation}\label{5.10}
0\leq \underline{\lim}_{t\to\infty}x(t) \leq \overline{\lim}_{t
\to\infty}x(t) \leq U \leq Y.
\end{equation}
Since condition $(C_k)$ holds for all $k\in J$, by lemma
\ref{Le5.5} there is a $\delta\in (0, 1]$ such that every solution
of $(\ref{1.1})$ in $\textup{int}{\R}^n_+$ satisfies
$\underline{\lim}_{t\to\infty}x(t) \geq \delta x^*$. Let
$\delta_0$ be the supremum of such $\delta$. Then, by lemma
\ref{Le5.2}, every solution of $(\ref{1.1})$ in
$\textup{int}{\R}^n_+$ satisfies
\[
\delta_0x^* = \hat{x} \leq \underline{\lim}_{t\to\infty}x(t) \leq
\overline{\lim}_{t\to\infty}x(t) \leq \hat{y},
\]
where $\hat{y}$ is given by $(\ref{5.1})$. By the assumption,
$x^*$ is on $\gamma_j$ for every $j\in I_n$. Hence, $\hat{x} = x^*
= \hat{y}$, so $x^*$ is a global attractor, if $\delta_0 = 1$.

Now suppose $\delta_0<1$. By lemma \ref{Le5.3}, every solution of
$(\ref{1.1})$ in $\textup{int}{\R}^n_+$ satisfies
\[
\hat{x} \leq \underline{\lim}_{t\to\infty}x(t) \leq
\overline{\lim}_{t\to\infty}x(t) \leq \hat{z},
\]
where $\hat{z}$ is given by $(\ref{5.2})$. For each $k\in J$,
since $\alpha_k\hat{x} = \delta_0\alpha_kx^* = \delta_0 <1$, if we
can show that
\begin{equation}\label{5.11}
\alpha_jy > 1, \quad j\in I_n\setminus\{k\}
\end{equation}
whenever $y\in [\hat{x}, \hat{x}^{\{k\}} + \hat{z}^{I_n
\setminus\{k\}}]\cap\gamma_k \not=\emptyset$, by lemma \ref{Le5.5}
again we have $\underline{\lim}_{t\to\infty}x(t) \geq \delta x^*$
for some $\delta\in (\delta_0, 1]$, which contradicts the
definition of $\delta_0$. This contradiction shows that $\delta_0
=1$.

To prove $(\ref{5.11})$ we introduce an affine mapping $f: [0, Y]
\to [0, Y]$ defined by $f(x) = \hat{x} + (1-\delta_0)x$. Then
$f(0) = \hat{x}$ and, for each $i\in I_n$, $f_i(Y) = \hat{x}_i +
(1-\delta_0)Y_i = \hat{x}_i +(1-\alpha_i\hat{x})a^{-1}_{ii} =
\hat{y}_i >\hat{x}_i$. Thus, $f([0, Y]) = [\hat{x}, \hat{y}]$ and,
for any $S\subset I_n$,
\[
f([0, Y^S]) = [\hat{x}, \hat{x}^{I_n\setminus S}+\hat{y}^S].
\]
Next we check that $f(U) = \hat{z}$ so that $F([0, U]) = [\hat{x},
\hat{z}]$ and, for any $S\subset I_n$,
\begin{equation}\label{5.12}
f([0, U^S]) = [\hat{x}, \hat{x}^{I_n\setminus S}+\hat{z}^S].
\end{equation}
For this purpose, we note that $x^*\in\cap_{i\in I_n}\gamma_i$ by
the assumption and $f(x^*) = \hat{x} + (1-\delta_0)x^* = x^*$. For
each $y\in [0, Y]\setminus\{x^*\}$, $f(y) = \delta_0x^*
+(1-\delta_0)y$ is on the line segment $\overline{yx^*}$. Thus,
$f(\gamma_i\cap[0, Y]) = \gamma_i\cap[\hat{x}, \hat{y}]$ for all
$i\in I_n$. For any fixed $i\in I_n$, if $U_i=a^{-1}_{ii}$ then,
by $(\ref{3.1})$, $a_{ij} =0$ or $Y^{\{i\}}$ is not below
$\gamma_j$, so $f(Y^{\{i\}}) = \hat{x} +(1-\delta_0)Y^{\{i\}} =
\hat{x}^{I_n\setminus\{i\}} +\hat{y}^{\{i\}}$ is not below
$\gamma_j$, for some $j\in I_n\setminus\{i\}$. By $(\ref{5.2})$ we
have $\hat{z}_i = \hat{y}_i$. Therefore, $f_i(U) = f_i(Y^{\{i\}})
= \hat{y}_i = \hat{z}_i$. If $U_i =0$, then $Y^{\{i\}}$ is below
$\gamma_j$ and every $y\in\gamma_i$ is on or below $\gamma_j$ for
all $j\in I_n\setminus\{i\}$. So
$\hat{x}^{I_n\setminus\{i\}}+\hat{y}^{\{i\}}$, which is on
$\gamma_i$, is below $\gamma_j$ and every $y\in
\gamma_i\cap[\hat{x}, \hat{y}]$ is on or below $\gamma_j$ for all
$j\in I_n\setminus\{i\}$. By $(\ref{5.2})$ again, $\hat{z}_i =
\hat{x}_i$. Hence, $f_i(U) = f_i(0) = \hat{x}_i = \hat{z}_i$. Now
suppose $0<U_i <a^{-1}_{ii}$. From $(\ref{3.1})$ we have $a_{ii} >
a_{ki}$ and $a_{ij} < a_{kj}$ for some $k, j\in
I_n\setminus\{i\}$. Then $\gamma_i\cap\gamma_k\cap[0, Y^{\{i,
j\}}]$ consists of a single point having
\[
\frac{a_{kj}-a_{ij}}{a_{ii}a_{kj}-a_{ij}a_{ki}}
\]
as its $i$th coordinate. Since $x\in \gamma_i\cap\gamma_k\cap[0,
Y^{\{i, j\}}]$ if and only if $f(x)=y\in
\gamma_i\cap\gamma_k\cap[\hat{x}, \hat{x}^{I_n\setminus\{i,
j\}}+\hat{y}^{\{i, j\}}]$ and $\alpha_i\hat{x} = \delta_0 =
\alpha_k\hat{x}$, by a routine check we see that
$\gamma_i\cap\gamma_k\cap[\hat{x}, \hat{x}^{I_n\setminus\{i, j\}}
+ \hat{y}^{\{i, j\}}]$ also consists of a single point $u^{(ijk)}$
with
\[
u_i^{(ijk)} = \hat{x}_i +
(1-\delta_0)\frac{a_{kj}-a_{ij}}{a_{ii}a_{kj} - a_{ij}a_{ki}}.
\]
From $(\ref{3.1})$ and $(\ref{5.2})$ we obtain
\begin{eqnarray*}
f_i(U) &=& \hat{x}_i + (1-\delta_0)\max\biggl\{\frac{a_{kj} -
a_{ij}}{a_{ii}a_{kj} - a_{ij}a_{ki}}: k, j \in I_n\setminus\{i\},
a_{kj} > a_{ij}\biggr\} \\
    &=& \max\{u_i^{(ijk)}: k, j \in I_n\setminus\{i\}, a_{kj} >
    a_{ij}\} \\
    &=& \hat{z}_i .
\end{eqnarray*}
Therefore, $f(U) = \hat{z}$. Now that $(\ref{5.12})$ is derived
from $f(U) = \hat{z}$, we must have
\[
f([0, U^{I_n\setminus\{k\}}]\cap\gamma_k) = [\hat{x},
\hat{x}^{\{k\}} +\hat{z}^{I_n\setminus\{k\}}]\cap\gamma_k.
\]
Then, for $y\in [\hat{x},
\hat{x}^{\{k\}}+\hat{z}^{I_n\setminus\{k\}}]\cap\gamma_k
\not=\emptyset$, there is a unique $x\in [0,
U^{I_n\setminus\{k\}}]\cap\gamma_k$ such that $f(x) = y$. By
condition $(C_k)$ for $k\in J$, $x$ is above $\gamma_j$ for all
$j\in I_n\setminus\{k\}$. Thus $y$ is above $\gamma_j$ for all
$j\in I_n\setminus\{k\}$, i.e. $y$ satisfies $(\ref{5.11})$.

If $Y$ replaces $U$ in the conditions and $\alpha_jx^* \geq 1$ is
assumed for all $j\in I_n \setminus J$, then $f(Y) \geq \hat{y}$
so $f([0, Y]) \supset [\hat{x}, \hat{y}]$. The above argument is
still valid after minor modifications.

Finally, that $x^*\leq U$ follows from $(\ref{5.10})$ and
$\lim_{t\to\infty}x(t) = x^*$ for every solution of $(\ref{1.1})$
in int${\R}^n_+$.

\vskip 5 pt \textit{Proof of theorem \ref{Th2.10}.} When $n\geq
m+2$, by assumption $(ii)$ and theorem \ref{Th2.1} with the
replacement of $I_n$ by $I_n\setminus J$ we know that
$(\ref{1.1})$ has a unique equilibrium $x^*\in \cap_{j\in J}\pi_j$
with $0< x^*_i \leq V_i$ for $i\in I_n\setminus J$. This is also
true with $x^* = Y^{\{k\}}$ when $I_n = J\cup\{k\}$ (i.e. $n =
m+1$). We next show that $x^*$ is a global attractor. From lemma
\ref{Le5.4} we know that every solution of $(\ref{1.1})$ in
int${\R}^n_+$ satisfies
\[
0=\hat{x} \leq \underline{\lim}_{t\to\infty}x(t) \leq
\overline{\lim}_{t\to\infty}x(t) \leq \hat{z}
\]
with $\hat{z} = V$. Then, by $(i)$ and $(ii)$, for each $k\in
I_n\setminus J$, either $V^{I_n\setminus\{k\}}$ is below
$\gamma_k$ or $\gamma_k\cap[0, V^{I_n\setminus\{k\}}]$ is above
$\gamma_j$ for every $j\in I_n\setminus\{k\}$. Applying lemma
\ref{Le5.5} to $(\ref{1.1})$ we obtain
\[
\lim_{t\to\infty}x(t)^J =0, \quad
\underline{\lim}_{t\to\infty}x_j(t)\geq \delta, \quad j\in
I_n\setminus J
\]
for some $\delta >0$ and all solutions of $(\ref{1.1})$ in
int${\R}^n_+$. Since $\{x^*\}$ is the largest invariant set of
$(\ref{1.1})$ in $\{y\in{\R}^n_+: y^J =0, y_j\geq\delta, j\in
I_n\setminus J\}$, $x^*$ must be a global attractor.

London Metropolitan University

Email: z.hou@londonmet.ac.uk

\end{document}